\documentclass[11pt]{amsart}
\usepackage{graphicx,amssymb,amsmath,amscd,verbatim}
\usepackage[all,cmtip]{xy}
\newtheorem{thm}{Theorem}[section]

\newtheorem{lem}[thm]{Lemma}
\newtheorem{prop}[thm]{Proposition}

\theoremstyle{definition}
\newtheorem{defn}[thm]{Definition}
\theoremstyle{remark}
\newtheorem{rem}[thm]{Remark}
\numberwithin{equation}{section}

\newcommand{\Real}{\mathbb R}

\newcommand{\Com}{\mathbb{C}}

\newcommand{\Kah}{K\"{a}hler }
\newcommand{\cxm}{\mathcal{M}}
\newcommand{\HB}{\mathcal{H}}

\newcommand{\dn}{d_\nabla}

\begin{document}

\title[Curvature for the complexified index cone of a cubic form]{A curvature formula for the complexified index cone of a cubic form}%
\author{Thomas Trenner}%
\address{Centre for Mathematical Sciences, Wilberforce Road, Cambridge CB3 0WA}%
\email{t.trenner@dpmms.cam.ac.uk}%


\date{\today}%
\begin{abstract}
We study the \Kah metric given by the logarithm of a
cubic form on its complexified index cone. Under mirror symmetry, this
metric should asymptotically correspond to the Weil--Petersson metric.
Using the theory of special \Kah manifolds, a proof of a curvature
formula for this metric is given.
\end{abstract}
\maketitle

\tableofcontents

\section{Introduction}

In this paper, we will prove a curvature formula for the Asymptotic Mirror Weil--Petersson metric --- or AMWP
metric --- on the complexified index cone of a cubic polynomial, as studied in \cite{TW}. 
This is motivated by classical mirror symmetry,
where on the complexified \Kah cone of a Calabi-Yau threefold we can define a metric corresponding, under the
mirror map, to the Weil--Petersson metric on the complex moduli space of a mirror manifold. To obtain the AMWP
metric, we neglect the contributions by instanton corrections. A curvature formula for the Weil--Petersson 
metric was given in \cite{Strominger1}, where it was also claimed that a similar result should hold on the
complexified \Kah cone. Further justification has been given in \cite{TW}, where the precise conjecture proved
in here was made (\cite[p.10]{TW}) and proved for the special case of ternary cubic polynomials. As the setting
in which Strominger proves the curvature formula for the Weil--Petersson metric in \cite{Strominger1} 
is special geometry, we will make use of the same setup here. \\

In the first two sections we give the basic definitions of special geometry and state the main results that
are needed to prove the curvature formula. In section 4, we quickly recall how special geometry is applied for
the complex moduli space of a Calabi--Yau threefold. The proof that the same curvature formula applies to the 
complexified index cone of a cubic polynomial will be given in section 5. This will include the examples of
complexified \Kah cones of the Calabi--Yau threefolds investigated in \cite{TW}.

\section{Affine special geometry}

In this section, we give the basic definitions of affine special \Kah manifolds and
derive a formula for the curvature tensor in terms of the prepotential. The main reference for this section is
Freed's paper on special \Kah manifolds \cite{Freed}, although we sometimes use slightly different notation and 
supply proofs where Freed omits them.

\begin{defn}
Let $(M, \omega)$ be a \Kah manifold. A \emph{special \Kah structure} on $M$ is a 
flat, torsion-free connection $\nabla$ on $TM$ such that
\begin{equation}\label{dnI} d_{\nabla}I = 0 \end{equation}
\begin{equation} \nabla \omega = 0 \end{equation}
Here $d_{\nabla}: \Omega^p(M;TM) \to \Omega^{p+1}(M;TM)$ 
denotes the associated exterior covariant derivative. 
\end{defn}

$\nabla$ and $\dn$ give the de Rham complex with coefficients in $TM$:
\begin{equation}\label{cxdR}
0\xrightarrow{}\Omega^0(TM)\xrightarrow{\dn}
\Omega^1(TM)\xrightarrow{\dn}
\Omega^2(TM)\xrightarrow{\dn}...
\end{equation}
Flatness of $\nabla$ is then expressed by requiring $\dn^2 = 0$. It implies that 
the Poincar\'e lemma is valid for (\ref{cxdR}).

\begin{lem}
We can express the torsion-free condition of $\nabla$ by $d_\nabla(\mathbf{id}_{TM}) = 0$, 
where we identify $\Gamma(End(TM))$ and $\Omega^1(TM)$.
\begin{proof}
Recall that the torsion tensor is defined as
\[ T(X,Y) = \nabla_X Y - \nabla_Y X - [X,Y]. \]
For the exterior covariant derivative we compute
\begin{eqnarray*} d_\nabla(\mathbf{id}_{TM})(X,Y) & = & \dn(\mathbf{id} (Y))X - \dn(\mathbf{id} (X))Y - 
\mathbf{id} ([X,Y]) \\
& = & \dn Y(X) - \dn X(Y) - [X,Y] \\
& = & \nabla_X Y - \nabla_Y X - [X,Y]
\end{eqnarray*}
\end{proof}
\end{lem}

Now choose a flat local framing $\{e^i\}, i=1,\ldots, n$, for $TM$ with dual framing
$\{\theta^i\}, i=1,\ldots, n$, for $T^*M$.

\begin{lem}
$d\theta^i = 0$.
\begin{proof}
In local coordinates, $id_{TM}$ is given by $\sum_i e^i \otimes \theta^i$. Then
$\dn(\mathbf{id}_{TM}) = \sum_i  e^i \otimes d\theta^i + \nabla e^i \otimes \theta^i$.
As the $e^i$ are a flat $TM$-framing, i.e. $\nabla e^i = 0$, the result follows.
\end{proof}
\end{lem}

\begin{defn}
We will call a coordinate system $\{ x_i, y_j\}, i,j=1,\ldots,n$ \emph{Darboux} if
\[ \omega = \sum_k dx_k \wedge dy_k . \]
\end{defn}

As our local coframing is exact, we have $\theta_i = dt_i$ for local coordinate functions
$t_i$. Because $\nabla \omega = 0$, we can choose the coordinate functions so that 
they form a Darboux coordinate system.
If we cover $M$ by local Darboux coordinate systems, the transition functions take the form
\[ \begin{pmatrix}x \\ y \end{pmatrix} = S \begin{pmatrix} \tilde{x} \\ \tilde{y} \end{pmatrix}
 + \begin{pmatrix} a \\ b \end{pmatrix}, \\ \ \ S \in Sp(2n;\Real), \ \ a,b \in \Real^n. \]

As giving a complex structure is equivalent to giving the decomposition of the complexified
tangent bundle in holomorphic and anti-holomorphic parts, instead of equation (\ref{dnI}) 
we can also require
\begin{equation} \dn \pi^{(1,0)} = 0, \end{equation}
where we identify $\pi^{(1,0)} \in Hom(TM^\Com, TM^{1,0})$ with its dual and use 
\[Hom(TM^{1,0}, TM^\Com) \cong \Omega^{1,0}(TM^\Com).\]

By the Poincar\'e lemma for (\ref{cxdR}) and $\dn \pi^{(1,0)} = 0$ there exists a complex 
vector field $\xi$ such that
\begin{equation}\label{vfxi} \nabla \xi = \pi^{(1,0)}. \end{equation}
In local Darboux coordinates,
\[ \xi = \frac{1}{2} \sum_{i=1}^n (z_i \frac{\partial}{\partial x_i} - w_i \frac{\partial}{\partial y_i}). \]
and thus
\[ \nabla \xi = \pi^{(1,0)} = \frac{1}{2} \sum_{i=1}^n(dz_i \otimes \frac{\partial}{\partial x_i} - dw_i \otimes \frac{\partial}{\partial y_i}). \]

This implies $Re(dz_i) = dx_i, Re(dw_i) = -dy_i$, which can be seen by taking the
real part of the last equation and using that the identity on $TM$ is the sum of the 
two projections: $\mathbf{id} = \pi^{1,0} + \pi^{0,1}$. Thus
\[ Re(\pi^{(1,0)}) = \frac{1}{2} \mathbf{id} = \frac{1}{2}(\sum_{i=1}^n(dx_i \otimes \frac{\partial}{\partial x_i} + dy_i \otimes \frac{\partial}{\partial y_i}). \]

\begin{defn}
Let $(M, \omega, \nabla)$ be a special \Kah manifold. A holomorphic coordinate
system $\{z_i\}$ is called \emph{special} if $\nabla Re(dz_i)=0$. Two special
coordinate systems $\{z_i\},\{w_j\}$ are called \emph{conjugate} if there exists
a flat Darboux system $\{x_i, y_j\}$ such that $Re(dz_i)=dx_i$, $Re(dw_i)= -dy_i$.
\end{defn}

By \cite[p.92]{ACD},this extension of the Darboux coordinate system to a
conjugate pair of special holomorphic coordinate systems is unique up to purely imaginary
translations.\\

The conjugate coordinate systems are related by a change of coordinates via 
holomorphic functions:
\begin{equation}
dw_i = \sum_j \tau_{ij} dz_j
\end{equation}

\begin{lem}\cite[p. 757]{Bartocci}
The holomorphic 2-form $\Theta = \sum_i dw_i \wedge dz_i$ is everywhere vanishing.
\end{lem}

This implies that $\theta = \sum_i w_i dz_i$ is closed, thus exact.

\begin{defn}
The holomorphic function $\mathcal{F}$ such that $\theta = d\mathcal{F}$ is
called the \emph{prepotential}.
\end{defn}

Let us choose the prepotential such that the following equations hold:
\begin{align}
\label{2.7}w_i &= \frac{4\partial \mathcal{F}}{\partial z_i} \\
\label{2.8}\tau_{ij} &= \frac{4\partial^2 \mathcal{F}}{\partial z_i \partial z_j}
\end{align} 

Then define $K = \frac{1}{2} Im(w_i \bar{z}_i) = 2 Im (\frac{\partial \mathcal{F}}{\partial z_i} \bar{z}_i)$.
We compute
\begin{align*}
 i\partial \bar{\partial}K &= 2i \partial \bar{\partial} Im (\frac{\partial \mathcal{F}}{\partial z_i} \bar{z}_i) \\
 &= 2i Im(\frac{\partial^2 \mathcal{F}}{\partial z_i \partial z_j})dz_i \wedge d\bar{z}_j\\
 &= -\frac{i}{2}Im(\tau_{ij})dz_i \wedge d\bar{z}_j 
\end{align*}

Note that this differs from the conventions in \cite{Freed}, where the \Kah potential is 
$\frac{1}{2} Im (\frac{\partial \mathcal{F}}{\partial z_i} \bar{z}_i)$.\\

On our \Kah manifold $M$ we have two affine connections: the Levi--Civita connection $D$ and
the flat, torsion-free connection $\nabla$. Define
\[ A_\Real = \nabla - D;\ \ \ \ \ \ A_\Real \in \Omega^1(M,End_\Real TM) \]
Extending $D$ and $\nabla$ to the complexified tangent bundle, we obtain $A$ and $\bar{A}$.

\begin{thm}\cite[p. 759, Thm. 3.9]{Bartocci}
\[ A \in \Omega^{1,0}(Hom(T^{1,0}, T^{0,1})) \]
\[ \bar{A} \in \Omega^{0,1}(Hom(T^{0,1}, T^{1,0})) \]
\end{thm}

For the curvature of the \Kah metric on $M$ we obtain the following formula:

\begin{thm}
$R_D = -(A \wedge \bar{A} + \bar{A} \wedge A)$
\begin{proof}
By definition, $\nabla = D + A_\Real$ is flat, so:
\begin{align}
 0 = \nabla^2 &=  (D + A_\Real)^2 = D^2 + D(A_\Real) + A_\Real(D) + A_\Real^2 \\
&= R_D + d_D(A) + A \wedge \bar{A} + \bar{A} \wedge A + A \wedge A + \bar{A} \wedge \bar{A} \\
&= R_D + d_D(A) + A \wedge \bar{A} + \bar{A} \wedge A ,
\intertext{as $A \wedge A$ and $\bar{A} \wedge \bar{A}$ vanish. We split the covariant exterior derivative of D: 
$d_D = \partial_D + \bar{\partial}_D$. As $D$ is the Levi--Civita connection corresponding to $\omega$,
we conclude that $\bar{\partial}_D = \bar{\partial}$. By decomposing the equation above into its 
$(r,s)$-components, we arrive at:}
0 &=  R_D + A \wedge \bar{A} + \bar{A} \wedge A + \partial_D(\bar{A}) + \bar{\partial}(A)
\intertext{We can further split this up into complex linear and anti-linear pieces, which both have to vanish:}
\label{affcurv}0 &=  R_D + A \wedge \bar{A} + \bar{A} \wedge A \\
0 &= \partial_D(\bar{A}) + \bar{\partial}(A)
\intertext{The complex linear part gives the equation we wanted. Also, in the last equation the 
two summands take values in different bundles, so we can conclude that $A$ is 
holomorphic and}
0 &= \partial_D(\bar{A}) 
\end{align}
Thus we arrive at the formula for $R$ we wanted.
\end{proof}
\end{thm}

\begin{defn}\label{cubicform}
We can also measure the failure of $\nabla$ to preserve the complex structure by 
the cubic form
\begin{equation}
\Xi = -\omega(\pi^{(1,0)}, \nabla \pi^{(1,0)}) = -\omega(\nabla \xi, \nabla^2 \xi),
\end{equation}
where we use (\ref{vfxi}) for the last equality, which of course only holds locally.
\end{defn}

\noindent In local coordinates we can then compute that the coefficients of $\Xi$ are given by the third 
partial derivatives of $\mathcal{F}$:

\begin{lem}
$\Xi_{ijk} = \frac{1}{2} \frac{\partial^3 \mathcal{F}}{\partial z_i \partial z_j \partial z_k}$
\begin{proof}
We may write $\pi^{(1,0)}$ in local coordinates as
\begin{align*}
\pi^{(1,0)} &= \sum_j \frac{\partial}{\partial z_j} \otimes dz_j. 
\intertext{Thus, using (\ref{2.7}) and (\ref{2.8}) we obtain}
\Xi &= -\omega \left(\sum_i dz_i \otimes \frac{\partial}{\partial z_i}, 
\nabla(\sum_j dz_j \otimes \frac{\partial}{\partial z_j})\right) \\
&= -\sum_{i,j,k,l,m}  \omega\left(\frac{1}{2}(\frac{\partial}{\partial x_i} - 
\tau_{im}\frac{\partial}{\partial y_m}),-\frac{1}{2} \frac{\partial \tau_{jl}}{\partial z_k}dz_k \otimes 
\frac{\partial}{\partial y_l}\right) \otimes dz_i \otimes dz_j\\\
&= \frac{1}{8} \sum_{i,j,k} \frac{\partial \tau_{ij}}{\partial z_k} dz_i \otimes dz_j \otimes dz_k \\
&= \frac{1}{2} \sum_{i,j,k} \frac{\partial^3 \mathcal{F}}{\partial z_i \partial z_j \partial z_k} dz_i \otimes 
dz_j \otimes dz_k.
\end{align*}

\end{proof}
\end{lem}

To compute the components of the curvature tensor in a local coordinate system $\{z_i\}_{i=1,\ldots,n}$, we use
the formula $0 =  R_D + A \wedge \bar{A} + \bar{A} \wedge A$ . We have:
\begin{align}
\label{omega}\omega &= \sum_{k,j} \frac{i}{2} g_{k\bar{j}} dz_k \wedge d\bar{z_j} \\
R_D &= \sum_{i,j,k,l} (R^j_{i k\bar{l}})_{i,j} dz_k \wedge d\bar{z_l} \\
(R_D)_{i\bar{j}k\bar{l}} &= \sum_m g_{m\bar{j}} R^m_{i k\bar{l}} \\
\Xi &= \sum_{i,j,k} \Xi_{ijk} dz_i \otimes dz_j \otimes dz_k \\
\intertext{As $A \in \Omega^{1,0}(Hom(T^{1,0}, T^{0,1}))$, we can expand it in a local
framing for the cotangent bundle:}
A &= \sum_i A_i dz_i, \ \ \ \ A_i \in Hom(T^{1,0},T^{0,1}) \\
\end{align}

We now show how, given a smooth cubic form $\Xi \in C^{\infty}(M,Sym^3 T^*M)$,
we recover the corresponding $\nabla$, or rather $A$.
In local coordinates $\omega$ is given by equation (\ref{omega}). 
Also, the curvature tensor $A$ is locally given by
\[ A = \sum_{i,j,k} A^{\bar{k}}_{ij} dz_i \otimes dz_j \otimes \frac{\partial}{\partial z_{\bar{k}}}. \]
Then
\begin{align*}
\Xi_{ijk} &= \frac{i}{2}\sum_{p,q} g_{p\bar{q}}dz_p \wedge d\bar{z}_q \left(\frac{\partial}{\partial z_i},
(A_\frac{\partial}{\partial z_j} \pi^{(1,0)})\frac{\partial}{\partial z_k} \right)\\
&= \frac{i}{2}\sum_{p,q} g_{p\bar{q}}dz_p \wedge d\bar{z}_q \left(\frac{\partial}{\partial z_i},
(\sum_{l,m} A^{\bar{m}}_{jl} dz_l \otimes \frac{\partial}{\partial \bar{z}_m})(\sum_j \frac{\partial}{\partial z_j} \otimes dz_j)\frac{\partial}{\partial z_k} \right)\\
&=  \frac{i}{2}\sum_{p,q} g_{p\bar{q}}dz_p \wedge d\bar{z}_q \left(\frac{\partial}{\partial z_i},
(\sum_{l,m} A^{\bar{m}}_{jl} dz_l \otimes \frac{\partial}{\partial \bar{z}_m})\frac{\partial}{\partial z_k} \right)\\
&=  \frac{i}{2}\sum_{p,q} g_{p\bar{q}}dz_p \wedge d\bar{z}_q \left(\frac{\partial}{\partial z_i},
\sum_{m} A^{\bar{m}}_{jk} \frac{\partial}{\partial \bar{z}_m}\right)\\
&= \frac{i}{2}\sum_q g_{i\bar{q}}A^{\bar{q}}_{kj}.
\intertext{Or equivalently}
A^{\bar{q}}_{jk} &= 2i\sum_i g^{\bar{q}i}\Xi_{ijk}.
\end{align*}

\noindent Using (\ref{affcurv}), we obtain the following formula for the curvature of $(M,\omega)$
in terms of the cubic form $\Xi$:

\begin{lem}
\begin{align}\label{affcurvfull}
 R_{i\bar{j}k\bar{l}} &= 4 \sum_{p,q} g^{p\bar{q}} \Xi_{ikp} \bar{\Xi}_{jlq} \\
 &= \sum_{p,q} g^{p\bar{q}} \mathcal{F}_{ikp} \bar{\mathcal{F}}_{jlq}
\end{align}
\end{lem}

\section{Projective special geometry}

Projective special \Kah manifolds are basically affine special \Kah manifolds
with a $\Com^*$-action. In Freed's setup this is achieved by taking a line bundle
$L$ over the \Kah manifold $M$. Then by deleting the zero section one obtains a
principal $\Com^*$-bundle, which naturally equips the total space with a 
$\Com^*$-action.

\begin{defn}
Let $(M, \omega)$ be a \Kah manifold, $L \to M$ a holomorphic line bundle. 
$(M, L, \omega)$ is called a \emph{Hodge manifold} if $c_1(L) = [ \omega ]$. This
implies that the \Kah class is integral.
\end{defn}

As $(L, \omega_L)$ is in particular a vector bundle, it has a global zero section. Removing this section,
we obtain a principal bundle $\tilde{M} \xrightarrow{\pi} M$ with structure group $\Com^*$. 
Every holomorphic line bundle admits a hermitian metric, which implies the existence of the \emph{Chern connection} $\nabla_L$ on $L$, 
which preserves the hermitian structure of the line bundle and has the property that
\[ \nabla_L^{0,1} = \bar{\partial}. \]
$\nabla_L$ induces a connection on the principal bundle $\tilde{M}$.\\

\begin{defn}
Given a vector bundle $\pi: V \to M$ and a continuous map $f: M' \to M$, we define the 
\emph{pullback bundle}
\[ f^*V := \{  (p,v) \in M' \times V | f(p) = \pi(v) \} \subseteq M' \times V \]
together with the natural projection to the first factor: $\pi': f^*V \to M'$.
For a section $s \in \Gamma(M,V)$ we define a \emph{pullback section} 
$f^*s \in \Gamma(M',V)$ by $f^*s = s \circ f$.
\end{defn}

The pullback bundle $\pi^*L \to \tilde{M}$ has a canonical nonzero holomorphic 
section $s$.

\noindent We define a pseudo-\Kah metric on $\tilde{M}$:

\begin{defn}
Let $\tilde{\omega} \in \Omega^{1,1}(\tilde{M})$ denote the form which is given by:
\begin{enumerate}
\item $\tilde{\omega}(U,X) = 0$
\item $\tilde{\omega}(X,Y) = |s|^2 \pi^*\omega(X,Y)$
\item $\tilde{\omega}(U,V) = -\frac{i}{2\pi}\partial \bar{\partial} |s|^2(U,V)$,
\end{enumerate}
where $X,Y$ are horizontal vector fields and $U,V$ are vertical vector fields.
\end{defn}

In the following we will establish a curvature formula for $\tilde{\omega}$ in terms of
$\omega$ and the cubic form $\Xi$ on $M$.

\begin{lem}
$(\tilde{M}, \tilde{\omega})$ is a pseudo-\Kah manifold. Furthermore
\[ \tilde{\omega} = \frac{i}{2\pi}\partial\bar{\partial}|s|^2 \]
and
\[ \pi^* \omega = \frac{i}{2\pi}\partial\bar{\partial}\log|s|^2 \]
\begin{proof}

The horizontal lift of a vector field $X \in TU$ is $X - \lambda h^{-1}\partial h(X) \frac{\partial}{\partial \lambda}$. Also, as by this formula $\tilde{\omega}$ is exact and real of type $(1,1)$, it defines a pseudo-\Kah metric of
signature $(\dim(M),1)$.
on $\tilde{M}$. It will not define a \Kah metric as it will be negative on vertical vector fields.
The second formula is just the standard curvature formula for the Chern connection of a holomorphic line
bundle, up to the constant factor of $\frac{i}{2\pi}$. 
\end{proof}
\end{lem}

Using this characterisation of the metrics, we can compute the \Kah metric
$\tilde{g}$ on $\tilde{M}$ in terms of the \Kah metric $g$ on $M$. We choose 
local coordinates $z_i$ on $U \subset M$ and a coordinate $\lambda$ on $L|_U$.
Choosing a nonzero local section $t$ of $L|_U$, we get local coordinates on
$\pi^*(U) \subset \tilde{M}$ by $(z,\lambda) \mapsto \lambda t(z), \lambda \in 
\Com^*$. 

\begin{lem}
Let $g, \tilde{g}$ be the \Kah metric on $M$ and $\tilde{M}$. Let $K_i :=
\partial_i \log K$, $K$ the \Kah potential of $\tilde{\omega}$. Also, we use the 
index $0$ in the metric for the fibre coordinate. Then
\begin{align}
 &\tilde{g}_{i\bar{j}} = K(-g_{i\bar{j}} + K_i K_{\bar{j}}) \\
 &\tilde{g}_{0\bar{i}} = K\lambda^{-1} K_i \\
 &\tilde{g}_{0\bar{0}} = (\lambda \bar{\lambda})^{-1}K \\
 &\tilde{g}^{i\bar{j}} = -K^{-1}g^{i\bar{j}} \\
 &\tilde{g}^{0\bar{i}} = \bar{\lambda}K^{-1}\sum_{k=1}^n g^{k\bar{i}}K_{\bar{k}} \\
 &\tilde{g}^{0\bar{0}} = K^{-1}\lambda \bar{\lambda}(1-\sum_{i,j=1}^n K_i K_{\bar{j}} g^{\bar{j}i})
\end{align}
\end{lem}

The entries of the inverse metric have been found by fixing the $\tilde{g}^{i\bar{j}}$ for $i,j = 1,\ldots,n$
to give the right answer when computing the Strominger formula in the projective case. This determines the other
entries and indeed gives an inverse matrix to $\tilde{g}$.\\ 

A tedious but easy calculation leads to expressions for the connection coefficients
of the Levi--Civita connection of $\tilde{g}$ in terms of the coefficients of the 
Levi--Civita connection for $g$. As both metrics are \Kah, we only have to compute
the ones which might be non-vanishing:

\begin{lem}
With notation as in the previous lemma, we compute the connection coefficients 
$\tilde{\Gamma}^i_{jk}$, using $\Gamma^i_{jk} = g^{i\bar{l}}\partial_j g_{\bar{l}k}$:
\begin{align}
&\tilde{\Gamma}^i_{jk} = \Gamma^i_{jk} + K_j \delta^i_k + K_k \delta^i_j \\
&\tilde{\Gamma}^i_{j0} = \lambda^{-1}\delta^i_j \\
&\tilde{\Gamma}^0_{ij} = \lambda \Gamma^k_{ij}K_k + 2\lambda K_i K_j - K^{-1}\lambda \partial_i \partial_j K\\
&\tilde{\Gamma}^0_{00} = \tilde{\Gamma}^0_{i0} = \tilde{\Gamma}^i_{00} = 0 
\intertext{It immediately follows that}
&\Gamma^i_{jk} = \tilde{\Gamma}^i_{jk} - \lambda \tilde{\Gamma}^i_{j0}K_k - \lambda \tilde{\Gamma}^i_{0k}K_i
\end{align}
\end{lem}

We can now compute the $(3,1)$ curvature tensor $R^m_{ik\bar{l}}$ of the 
Levi--Civita connection on $M$, using the following formula:
\begin{align*}
R^m_{ik\bar{l}} &= -\partial_{\bar{l}} \Gamma^{m}_{ik} 
\intertext{Expressed in terms of the $(3,1)$ curvature tensor 
$\tilde{R}^m_{ik\bar{l}}$ on $\tilde{M}$ and the \Kah potential $K$, we obtain:}
R^m_{ik\bar{l}} &= -(\tilde{R}^m_{ik\bar{l}} - \partial_{\bar{l}}(\delta^m_i 
K_k + \delta^m_k K_i)),
\intertext{where $i,j,k,m \in \{1,\ldots,n\}$. Expressing $\tilde{R}^m_{ik\bar{l}}$ in terms of the
cubic form, we obtain:}
R^m_{ik\bar{l}} &= -4\sum_{p,q,m,n} \tilde{g}^{p\bar{q}} \tilde{g}^{n\bar{m}} \tilde{\Xi}_{ikp} \tilde{\bar{\Xi}}_{mlq} + \partial_{\bar{l}}(\delta^m_i K_k + \delta^m_k K_i)
\end{align*}
\begin{align*}
\intertext{Plugging in (\ref{affcurvfull}), 
we obtain for the $(4,0)$ curvature tensor:}
R_{i\bar{j}k\bar{l}} &= \sum_m g_{m\bar{j}}R^m_{ik\bar{l}} \\ 
&= \sum_m g_{m\bar{j}} (-\tilde{R}^m_{ik\bar{l}} +
\partial_{\bar{l}}(\delta^m_i K_k + \delta^m_k K_i)) \\
&= -\sum_m g_{m\bar{j}}4 \sum_{p,q,n} \frac{1}{K} g^{p\bar{q}} \frac{1}{K} g^{n\bar{m}} \tilde{\Xi}_{ikp} \tilde{\bar{\Xi}}_{mlq} + \sum_m g_{m\bar{j}}\partial_{\bar{l}}\delta^m_i K_k +
\sum_m g_{m\bar{j}} \partial_{\bar{l}} \delta^m_k K_i \\
&= -\frac{4}{K^2} \sum_{p,q} g^{p\bar{q}} \Xi_{ikp} \bar{\Xi}_{jlq} + g_{i\bar{j}}\partial_{\bar{l}} K_k +
g_{k\bar{j}} \partial_{\bar{l}} K_i \\
&= -\frac{4}{K^2} \sum_{p,q} g^{p\bar{q}} \Xi_{ikp} \bar{\Xi}_{jlq} + g_{i\bar{j}}g_{k\bar{l}} +
g_{k\bar{j}}g_{i\bar{l}} \\
&= -\frac{1}{K^2}\sum_{p,q} g^{p\bar{q}} \mathcal{F}_{ikp} \bar{\mathcal{F}}_{jlq} + g_{i\bar{j}}g_{k\bar{l}} +
g_{k\bar{j}}g_{i\bar{l}} 
\end{align*}

For the last equality we take the above description of $\tilde{\omega}$. Using this in definition
\ref{cubicform} and remembering that $\tilde{\nabla}$ is just the pullback of $\nabla$ by the $\Com^*$-action,
we see that $\tilde{\Xi}$ and $\Xi$ do not differ for the vector fields coming from $M$. It should be noted 
that for a vector field $X$ on $M$ we have to use the horizontal lift to $\tilde{M}$ to make this 
identification. The final step is to rewrite the \Kah potential in terms of the cubic form:

\begin{prop}
Denote the coordinates on $\tilde{M}$ by $\lambda, z_1,\ldots,z_n$ with\\ $z_j = x_j + i y_j$. In special coordinates $t_j = z_j/\lambda, j=1,\ldots ,n$, writing $\partial_j$ for $\frac{\partial}{\partial t_j}$, the K\"{a}hler potential for the Weil--Petersson metric is given by
\[ K= - \log i(\sum_j ((t_j - \bar{t_j})(\partial_j \mathcal{F'} + \bar{\partial_j}\bar{\mathcal{F'}})) + 2 \bar{\mathcal{F'}} - 2\mathcal{F'}) = 8 \mathcal{F}(y_1,\ldots,y_n)\]
\begin{proof}
Let the $t_j = \frac{z_j}{\lambda}$ form a local coordinate system and define $\mathcal{F'}(t) = \lambda^{-2}\mathcal{F}(z)$. Then observe that $\lambda \partial_j \mathcal{F}' = \frac{\partial \mathcal{F}}{\partial z_j}$. Going to inhomogeneous coordinates, we have to account for the following term that has to be corrected in the sum:
\[ \lambda \frac{\partial \mathcal{F}}{\partial \lambda} = \lambda^2 (2\mathcal{F'} -  \sum_j t_j \partial_j \mathcal{F'}) \]
Plugging this into the formula given in the Lemma above we obtain the result.
\end{proof}
\end{prop}

\begin{defn}\cite[p.46]{Freed}
Let $(M, L, \omega)$ be an $n$-dimensional Hodge manifold. A \emph{projective
special \Kah structure} is a triple $(V, \nabla, Q)$, where
\begin{enumerate}
\item $V \to M$ is a rank $n+1$ holomorphic vector bundle with a holomorphic 
inclusion $L \hookrightarrow V$.
\item $\nabla$ is a flat connection on $V_\Real$, the underlying real vector bundle, 
such that the extension of $\nabla$ to the complexification $(V_\Real)_\Com$ 
satisfies $\nabla(L) \subset V$ and the section
\begin{align*}
M &\to \mathbb{P}[(V_\Real)_\Com] \\
m &\mapsto L_m 
\end{align*}
is an immersion with respect to $\nabla$, or equivalently, $L$, seen as a section of 
$\mathbb{P}[(V_\Real)_\Com]$, is transverse to the horizontal distribution of the connection.
\item $Q$ is a non-degenerate skew-symmetric form on $V_\Real$ of type $(1,1)$
with respect to the complex structure. Moreover, it should be flat with respect to 
$\nabla$. Finally, taking the extension to $(V_\Real)_\Com$, we assume 
$Q|_{L \times \bar{L}} = \frac{i}{2\pi}\omega_L$.
\end{enumerate}
\end{defn}

\begin{rem}
By $\nabla(L) \subset V$ we mean that if we take a section $s \in \Gamma(L)$, where L
is embedded in $V$, we get $\nabla(s) \in \Gamma(V) \otimes \Omega_M$.
\end{rem}

A projective special \Kah structure defines a certain kind of variation of Hodge structures. 
If we furthermore require the existence of a lattice in $V_\Real$ such that $Q$ restricted
to the lattice takes integer values, the corresponding variation of Hodge structures will be
polarised.\\

\noindent The connection to affine special geometry is established by the following proposition:

\begin{prop}\cite[p.46]{Freed}
Let $(M, L, \omega)$ be a Hodge manifold with associated pseudo-\Kah manifold 
$(\tilde{M}, \tilde{\omega})$ and canonical holomorphic non-zero section $s$.
Then a projective special \Kah structure on $(M, L, \omega)$ is equivalent to a
$\Com^*$-invariant special pseudo-\Kah structure $\tilde{\nabla}$ on 
$(\tilde{M}, \tilde{\omega})$ with $\tilde{\nabla}s = \pi^{(1,0)}$.
\end{prop}

\noindent The condition $\tilde{\nabla}s = \pi^{(1,0)}$ is needed to guarantee that the map in part (2) of
the definition of a projective \Kah structure is an immersion: consider $s$ as a holomorphic vertical vector field on 
$\tilde{M}$ and denote by $\tilde{\nabla}$ the pullback of the connection to $\pi^* V$. Then the map
\begin{align}\label{connmap}
\tilde{\nabla} s: T\tilde{M} &\to \pi^* V \\
\nonumber \partial / \partial z_i &\mapsto \nabla_{ \partial / \partial z_i} s
\end{align}
is an isomorphism if and only if $s$ is transverse to the horizontal distribution of the connection. The
underlying real isomorphism induces a real flat connection on $ T\tilde{M}$.\\

If we denote the sections of $V$ by $s,l_1,\ldots,l_n$, the immersion condition says that the effect of the 
connection on $s$ is given by:
\[	\tilde{\nabla} s = \sum_{j=0}^{n} \frac{\partial}{\partial z_j} \otimes dz_j \]
This implies that under (\ref{connmap}), $s$ corresponds to a complex vector $\zeta$ field that satisfies 
$\tilde{\nabla} \zeta = \pi^{(1,0)}$ because of $\pi^{(1,0)} = \sum_j \frac{\partial}{\partial z_j} \otimes 
dz_j.$\\

Another useful fact is that $\tilde{Q}$, the pullback of the skew-symmetric form via $\pi$, pulls back to 
$-\tilde{\omega}$ under the isomorphism, which can be seen by differentiating the equation we assume in point
(iii) of the projective special \Kah structure.
Using this together with definition \ref{cubicform} and the fact that $\omega(\zeta,\tilde{\nabla}\zeta) = 0$ 
as \(\zeta\) is a holomorphic $(1,0)$ vector field, we conclude that the cubic form is given by
\begin{equation}\label{projcube}
\Xi = -\omega(\zeta,\tilde{\nabla}^3\zeta) = \tilde{Q}(\tilde{\nabla}^3 s, s). 
\end{equation}

In the following sections we will describe two instances of projective
special \Kah manifolds that arise in Mirror Symmetry, the complex moduli space and the (complexified) \Kah 
moduli space of a Calabi--Yau threefold. A curvature formula for the former can be found in \cite{Strominger1},
although it is derived in a slightly different fashion. We want to compute the curvature of the complexified
\Kah moduli space of a Calabi--Yau threefold with the line bundle $\mathcal{H}^{0,0}$. We do that by proving the
formula for the more general case of a complexified index cone of a cubic form.\\

\section{Complex moduli space as a Special \Kah Manifold}

Let $X$ be a Calabi--Yau threefold, $\cxm$ the moduli space of complex structures. Over $\cxm$,
consider the Hodge bundle $\HB$ with fibre $H^3(X,\Com)$ and Hodge filtration $\mathcal{F}^\bullet$. This gives a 
$2n+2$-dimensional vector bundle over $\cxm$. Using integral cohomology,
we can equip it with a flat connection, called the \emph{Gauss--Manin connection}.
There exists a symplectic basis for 
$H^3(X,\Com)$ and this provides local flat sections for $\cxm$ with respect 
to the Gauss--Manin connection on $\HB$. These bases
are unique up to $Sp(2n+2,\Real)$ transformations. Thus $\HB$ is a flat, 
holomorphic $Sp(2n+2,\Real)$ bundle.\\

As all the information about the Hodge structure is already contained in $H^{3,0} \oplus H^{2,1}$, we get an 
$n+1$-dimensional sub-bundle $\mathcal{V} \cong \mathcal{F}^2$. 
On $\HB$, we have two metrics: One is a hermitian metric given by the polarisation of Hodge structure 
$Q(C\cdot,\cdot)$ 
(in most accounts called Hodge metric), where $C$ is the Weil operator, i.e. 
acting by multiplication by $i^{p-(n-p)}$ on $H^{p,q}$. The other one is given just by taking 
the polarisation with a leading factor of $-i^n$, which in general is not positive 
definite. We will call it the \emph{indefinite metric}. These two metrics differ 
by $(-1)^p$ on $H^{p,q}$. While the indefinite metric is by definition preserved 
by the Gauss--Manin-connection, the Hodge metric is not. This implies that $\nabla^{GM}$
is not the Chern connection of the hermitian line bundle $\mathcal{H}^{3,0}$.\\

On any Calabi--Yau manifold of the family described by the complex moduli space, the nowhere-vanishing 
holomorphic 3-form $\Omega$ is only defined up to multiplication by a non-zero complex number. So in this
family, we can multiply $\Omega$ by a non-vanishing holomorphic function $f(z)$. As we want no-where 
vanishing 3-forms, we take the following transformation for $\Omega$:
\[ \Omega' = e^f(z) \Omega, \]
where $f$ is a continuous function of the moduli space coordinates. Now at different points in the complex
moduli space, we have different complex structures and thus a different decomposition of the tangent bundle of
our Calabi--Yau threefold in holomorphic and anti-holomorphic parts. This means that $(3,0)$-forms that
are holomorphic in one complex structure may not stay holomorphic in another complex 
structure. As $\Omega$ is by definition a holomorphic form, it defines us a section of 
$H^3$. The multiples of $\Omega$ define the line bundle $\mathcal{H}^{3,0}$
which is a sub-bundle of both $\HB$ and $\mathcal{V}$. 
We also note for later that because of the freedom in 
choice for \(\Omega\) we obtain a section of the projectivization of $\mathcal{V}$.\\

\noindent We equip $\cxm$ with a projective special \Kah structure $(\mathcal{V},\nabla,Q)$ given by:
\begin{enumerate}
\item 
The vector bundle $\mathcal{V} := \mathcal{H}^{3,0} \oplus \mathcal{H}^{2,1}$ 
with the holomorphic inclusion  $L := \mathcal{H}^{3,0} \hookrightarrow \mathcal{H}^{3,0} 
\oplus \mathcal{H}^{2,1}$.
\item 
The Gauss--Manin connection $\nabla^{GM}$ on the underlying real bundle $\mathcal{V_\Real}$. For the rest of 
this chapter, we will denote the Gauss--Manin connection simply by $\nabla$. 
\item 
The non-degenerate form $Q$ is the polarisation of the variation of Hodge  structures on the complex moduli 
space, divided by $2\pi$. Then $Q$ is flat with respect to $\nabla^{GM}$.
\end{enumerate}

\begin{lem}
$(\mathcal{V},\nabla,Q)$ defines a projective special \Kah structure.
\begin{proof}
The holomorphicity of the inclusion $L \hookrightarrow \mathcal{V}$ follows from a theorem of Griffiths 
\cite[p.24]{Griffiths1} stating that the bundles in the Hodge filtration $\mathcal{F}^p$ corresponding to $\HB$ 
are holomorphic sub-bundles of the full Hodge bundle $\HB$). Here $\mathcal{F}^p = \sum_{i=p}^{3} \mathcal{H}^{p,3-p}$.
The requirement that $\nabla(L) \subset \mathcal{V}$ after complexification of the underlying real bundles is just 
Griffiths transversality for the Hodge bundles.  For the immersion condition in part (2) of the definition of a projective special \Kah structure,
we note that $\mathcal{V}$ is generated by $\Omega$ and $\nabla_{\frac{\partial}{\partial z_i}}\Omega$, 
$i=1,\ldots,n$. But $L$ is generated by the non-vanishing $(3,0)$-form $\Omega$ and 
$\nabla_{\frac{\partial}{\partial z_i}}\Omega$ is a non-zero section of $\mathcal{V}$, so $L$, seen as a section of 
$\mathbb{P}[(\mathcal{V}_\Real)_\Com]$, is transverse to the horizontal distribution of the connection.
The remaining required properties for $Q$ are immediate from the definition.
\end{proof}
\end{lem}

With this setup and the calculations of the previous section, we obtain the following curvature formula for the
complex moduli space appearing in Strominger's paper \cite{Strominger1}:
\[ R_{i\bar{j}k\bar{l}} = g_{i\bar{j}}g_{k\bar{l}} +
g_{i\bar{l}}g_{k\bar{j}} - e^{2K} \sum_{p,q} g^{p\bar{q}}
Y_{ikp}\overline{Y_{jlq}}, \]
where $g$ is the Weil--Petersson metric, given by taking as \Kah potential $\log Q(\Omega, \bar{\Omega})$, and $Y_{ijk}$ denotes the 
Yukawa couplings (see \cite[p.102]{CoxKatz1} for details). We use (\ref{projcube}) to obtain as the cubic form 
$\Xi$ the well-known formula for the Yukawa couplings. The cubic form does not fix the prepotential completely, 
but as the curvature only depends on $\Xi$, we obtain the well-known original formula for the complex moduli 
space.\\

Proofs of this formula for the Weil--Petersson metric on the complex moduli space of a Calabi--Yau threefold have 
appeared in the literature before. A proof relying on the earlier extrinsic definition of special geometry via 
special coordinates and their behaviour under coordinate changes can be found in \cite{Strominger1}. Two more 
proofs, using different methods --- most importantly a reformulation of the problem in terms of the period 
mapping and a theorem by Griffiths on the curvature of Hodge bundles --- have been found by Wang \cite{Wang}.

\section{Complexified index cone of a cubic form as a projective Special \Kah manifold}

Take any cubic polynomial with real coefficients $f(y_1,\ldots,y_n)$ defined on $\Real^n$. 

\begin{defn}
The \emph{index cone} of $f$ is the open cone $W \subset \Real^n$ where f is positive and the
Hessian matrix $(\partial^2 f / \partial y_i \partial y_j)$ has index $(1,n-1)$.
\end{defn}

\begin{defn}
The \emph{complexified index cone} $M$ is given by $(\Real^n + iW)/im(\mathbb{Z}^n)$, where $\mathbb{Z}^n$ is 
mapped into the first summand.
\end{defn}

Let $t_i = x_i + i y_i$. As was shown in \cite[p.7ff.]{TW}, starting with the complexification of a cubic 
polynomial f --- e.g. the cubic intersection form of a Calabi--Yau threefold --- we get a \Kah potential from it by
using the formula
\[ K = -\log i(\sum_j (t_j - \bar{t_j})(\partial_j f + \bar{\partial_j} \bar{f}) + 2\bar{f} - 2f). \]
This potential function is independent of the $x_i$, and is just a multiple of the original polynomial $f$.

We consider the Hessian metric $(- \frac{1}{4} \partial^2 \log f / \partial y_i \partial y_j)$ on the complexified index cone,
obtained by taking $K$ as the \Kah potential. By \cite[Lemma 2.1]{TW}, this defines a \Kah metric on the 
whole complexified index cone.

It was conjectured in \cite[p.10]{TW} that the following equivalent of the Strominger formula holds for the 
Hessian metric on the complexified index cone:
\begin{equation}\label{Strom2}
R_{i\bar{j}k\bar{l}} = g_{i\bar{j}}g_{k\bar{l}} +
g_{i\bar{l}}g_{k\bar{j}} -  \sum_{p,q} g^{p\bar{q}}\frac{f_{ikp}f_{jlq}}{64 f^2}
\end{equation}

The factor of 64 is due to the fact that $K = 8\mathcal{F}(y) = 8f$, so $\frac{1}{K^2} = \frac{1}{64f^2}$.\\

In the following we will show that the complexified index cone together with the Hessian metric given by the 
cubic polynomial $f$ and the trivial complex line bundle $L=\Com \times M$ form a projective special \Kah 
manifold. Denote the total space obtained by deleting the zero section from $L$ again by $\tilde{M}$. We choose 
local coordinates $z_1,\ldots,z_n,\lambda$, where $\lambda$ denotes the fibre coordinate.

\noindent Then a projective special \Kah structure $(\mathcal{V},\nabla,Q)$ on $\tilde{M}$ given by:
\begin{enumerate}
\item 
The vector bundle $\mathcal{V} := \Com^{n+1}$ 
with the obvious holomorphic inclusion  of $L = \Com$. The underlying real bundle is just $\Real^{n+1}$.
\item 
Choose rational non-vanishing sections $s,s_1,\ldots,s_n$ for $V$, with $s$ the section for $L$. Our aim is to 
identify these sections with the basic vector fields $\partial / \partial z_i, \partial / \partial \lambda$, 
where the former should be identified with  $s_1$ up to $s_n$ and the last one should be identified with $s$. 
We use the shorthand notation $\partial z$ for the vector field $\partial / \partial z$.
On the real bundle define a connection by setting
\begin{align*}
\nabla_{\partial z_i}s&= s_i\\
\nabla_{\partial z_i}s_j&= 
\sum_{k=1}^n \frac{\partial^3 f}{\partial z_i \partial z_j \partial z_k}\bar{s}_k\\
\nabla_{\partial z_i}\bar{s}_j&= \delta_j^i \bar{s}\\
\nabla_{\partial z_i}\bar{s}&= 0,
\end{align*}
where $\bar{s_k}$ are rational sections of $(\Com^n)^\vee$, the dual to $V\backslash L$.
\item 
The non-degenerate form $Q$ is given by $Q(\alpha,\beta) = \frac{1}{2} \alpha \beta$. 
\end{enumerate}

We compute 
\begin{align*}
\Xi_{ijk} &= Q(\tilde{\nabla}_{\partial z_i} \tilde{\nabla}_{\partial z_j} 
\tilde{\nabla}_{\partial z_k} s, s)\\ &= Q(\tilde{\nabla}_{\partial z_i} \tilde{\nabla}_{\partial z_j} 
s_k, s)\\ &= Q(\tilde{\nabla}_{\partial z_i} \sum_{i=1}^n 
\frac{\partial^3 f}{\partial z_i \partial z_j \partial z_k} \bar{s}_i, s)\\ 
&= Q(\frac{\partial^3 f}{\partial z_i \partial z_j \partial z_k} \bar{s}, s)\\
&= \frac{1}{2}\frac{\partial^3 f}{\partial z_i \partial z_j \partial z_k} 
\end{align*}
which shows that (\ref{Strom2}) holds.

\bibliographystyle{amsplain}

\end{document}